\documentclass[submission%
,nohyperref
]{dmtcs}
\usepackage[latin1]{inputenc}
\usepackage[round]{natbib}
\usepackage{amsmath,tikz}
\let\set\mathbb
\def\sgn{\operatorname{sgn}}

\newtheorem{thm}{Theorem}

\long\def\family#1#2#3{%
  Family #1\rule{0pt}{1.1em}\hfill\vspace{5pt}\break
  \begin{tabular}{@{}p{.79\hsize}@{}p{.21\hsize}@{}}
    Defining equations:\hfill\break
    #2
    &
    Example:\hfill\break
    #3
  \end{tabular}\vspace{2pt}%
}

\author{Manuel Kauers\addressmark{1}\thanks{Email: \email{mkauers@risc.jku.at}. Partially supported by the Austria FWF grants Y464-N18 and F50-04.}
  \and Rika Yatchak\addressmark{1}\thanks{Email: \email{ryatchak@risc.jku.at}. Partially supported by the Austria FWF grant F50-04}}
\title{Walks in the Quarter~Plane with Multiple~Steps}
\address{\addressmark{1}RISC, Johannes Kepler University, Linz, Austria}
\received{tba}
\revised{tba}
\accepted{tba}
\begin{document}
\maketitle

\begin{abstract}
  We extend the classification of nearest neighbour walks in the quarter plane
  to models in which multiplicities are attached to each direction in the step
  set. Our study leads to a small number of infinite families that completely
  characterize all the models whose associated group is D4, D6, or~D8. These
  families cover all the models with multiplicites 0, 1, 2, or~3, which were 
  experimentally found to be D-finite --- with three noteworthy exceptions. 
\end{abstract}

\section{Introduction}

We consider quadrant walk models where step sets may contain
several distinguishable steps pointing into the same direction. For example, the
step sets $\{\leftarrow,\downarrow,\nearrow\}$ and
$\{\leftarrow,\leftarrow',\downarrow,\nearrow\}$ are considered different, as
the latter contains two different ways of going to the left. The objects being
counted are then walks in the quarter plane starting at the origin, consisting
of $n$ consecutive steps taken from the step set in such a way that the walk
never leaves the first quadrant, ending at a point $(i,j)\in\set N^2$, and one
of $k$ different colors is attached to each step in the walk whose multiplicity
in the step set is~$k$. For each model (viz., for each multiset of admissible
directions), we want to know whether the corresponding generating function
$f(x,y,t)=\sum_{n=0}^\infty \sum_{i,j} f_{i,j,n} x^i y^j t^n$ which counts the
number $f_{i,j,n}$ of walks of length~$n$ ending at $(i,j)$ is D-finite. As
usual, a power series in $t$ is D-finite if it satisfies an ordinary linear
differential equation with polynomial coefficients.

If we let $a_{u,v}$ denote the multiplicity of the direction
$(u,v)\in\{-1,0,1\}^2\setminus\{(0,0)\}$, then the number $f_{i,j,n}$ of walks
of length~$n$ ending at $(i,j)$ is uniquely determined by the recurrence
equation
\[
   f_{i,j,n+1} = \sum_{u,v} a_{u,v} f_{i-u,j-v,n}\qquad(n\in\set N, i,j\in\set N)
\]
together with the initial values $f_{0,0,0}=1$, $f_{i,j,0}=0$ for $(i,j)\neq(0,0)$,
and the boundary conditions $f_{-1,j,n}=f_{i,-1,n}=0$ for all $i,j,n$. Equivalently,
we can say that the generating function $f(x,y,t)=\sum_{n=0}^\infty \sum_{i,j} f_{i,j,n} x^i y^j t^n
\in\set Q[x,y][[t]]$ satisfies the functional equation
\begin{alignat}1\label{eq:funeq}
  &\Bigl(1 - t \sum_{u,v} a_{u,v}x^uy^v\Bigr)\,f(x,y,t)\notag\\
  &\quad{}= 1 - \frac1{ty}\Bigl(\sum_u a_{u,-1}x^u\Bigr) f(x,0,t) - \frac1{tx}\Bigl(\sum_v a_{-1,v}y^v\Bigr) f(0,y,t) + \frac{a_{-1,-1}}{txy} f(0,0,t).
\end{alignat}
Its first terms are
\begin{alignat*}1
 &f(x,y,t)= 1+ \bigl(a_{1,1}\,x y+ a_{1,0}\,x+ a_{0,1}\,y\bigr)t 
   +  \bigl(a_{1,1}^2\,x^2 y^2+2 a_{1,0} a_{1,1}\, x^2 y+(a_{1,0}^2+a_{1,-1} a_{1,1})x^2 \\
 &\qquad{}+2 a_{0,1} a_{1,1}\,x y^2+2a_{0,1} a_{1,0}\, x y
   + (a_{0,1} a_{1,-1}+a_{0,-1} a_{1,1})x
   +(a_{0,1}^2+a_{-1,1} a_{1,1})y^2\\
 &\qquad{}+(a_{-1,1} a_{1,0}+a_{-1,0} a_{1,1})y 
   +(a_{0,-1}a_{0,1}+a_{-1,0} a_{1,0}+a_{-1,-1} a_{1,1})\bigr)t^2 + \cdots. 
\end{alignat*}
This means, for example, that there are $a_{-1,1}a_{1,0}+a_{-1,0}a_{1,1}$ many walks
of length $n=2$ ending at $(i,j)=(0,1)$.

For the models where all multiplicities $a_{u,v}$ are in $\{0,1\}$, a complete
classification is available: among the $2^8=256$ different models, \cite{bousquet10} 
identified 79 nontrivial cases. For 22 of them
they prove that the generating function is D-finite using certain symmetry
groups~$G$ associated to each of the models. For a 23rd model, the notorious
Gessel model $\{\leftarrow,\to,\nearrow,\swarrow\}$, their techniques do not
apply but a proof by a different method based on computer algebra was found
by \cite{bostan10}. A computer-free proof was later found by \cite{bostan13}.  
The remaining 56 models are not D-finite: \cite{rechnitzer09} 
and \cite{melczer13} showed that the generating functions of 
five of these models have infinitely many singularities
and therefore are not D-finite. For the remaining models, \cite{bostan14a} proved that the counting sequences $f_{0,0,n}$ for walks returning
to the origin have asymptotic behaviour for $n\to\infty$ that D-finite
functions cannot possibly have.

\def\vec#1#2#3{\Bigl(\vcenter{\baselineskip=0pt\hbox{$\scriptstyle\mathstrut#1$}\kern-1pt
    \hbox{$\scriptstyle\mathstrut#2$}\kern-1pt\hbox{$\scriptstyle\mathstrut#3$}}\Bigr)}
The need for a classification of quarter plane models with multiplicities arose
in the classification project for octant models in 3D~\citep{bostan14}, as it turns
out that some models in 3D can be reduced by projection to 2D models with
multiplicities. For example, it is easy to see that the generating function for
the octant model with step set
$\{\vec{-1}{\;0}{\;0},\vec{-1}{\;0}{\;1},\vec{\;0}{-1\!}{\;0},\vec110\}$ is
D-finite if and only if the quadrant model with step set
$\{\leftarrow,\leftarrow',\downarrow,\nearrow\}$ is. \cite{bostan14}
have classified only the 527 models that they needed for their study, and point
out that the classification problem for models with multiplicities is of
interest in its own right. 

For the present paper we carried out a systematic search
over all the $4^8=65536$ models where each of the eight directions may have any
of the four multiplicities 0,~1,~2,~3. Of these, 30307 are nontrivial and
essentially different, and of these, 1457~turn out to be D-finite, and of these,
79~are even algebraic.  Going one step further, we have identified families of
D-finite models in which some or all of the ``multiplicities'' are arbitrary
complex numbers. Rather than asking for a fixed model what the corresponding group
is, we ask for a fixed group what all the models leading to this group are. In
this way we obtain a small number of families that completely characterize all
the models which lead to groups with at most eight elements. This
characterization covers 1454 of the 1457 D-finite cases we discovered for
multiplicities in $\{0,1,2,3\}$, the remaining three models have a group of
order~10, which was too hard for us to analyze in full generality. In view of the
fact that all models previously considered had either a finite group of order at most eight
or an infinite group, the appearance of these models was a surprise to us. We were
less surprised to find, after spending some 6.5~years of computation time, that none of the
models with a (probably) infinite group appears to be D-finite based on the inspection
of the first 5000 terms. 

For models with multiplicities $a_{u,v}\in\{0,1\}$ it is noteworthy that the generating
function for a model is D-finite if and only if the associated group is finite, and it
is algebraic if and only if the so-called orbit-sum (cf. Section~\ref{sec:d4} below) is zero.
It seems that these equivalences remain true for models with multiplicities. 


\section{Models of Interest}

Our reasoning largely follows that of \cite{bousquet10}. Their
first step is to identify the interesting models. By a model, we understand here a particular
choice of multiplicities $a_{u,v}\in\set C$ (not necessarily integers). For each such model
there is a corresponding generating function $f(x,y,t)\in\set C[x,y][[t]]$, and we want to
identify the models whose generating functions are D-finite. 

A model is uninteresting if $a_{1,-1}=a_{1,0}=a_{1,1}=0$ or 
$a_{-1,1}=a_{0,1}=a_{1,1}=0$ or 
$a_{-1,-1}=a_{-1,0}=a_{-1,1}=0$ or
$a_{-1,-1}=a_{0,-1}=a_{1,-1}=0$, because in either of these cases the corresponding generating
function is algebraic and it is well-understood why \citep[Section~VII.8]{flajolet09}. 
Secondly, if two models can be obtained from one another by reflecting the step set about the 
diagonal~$x=y$, then the corresponding generating functions can be obtained from one another
by exchanging the variables $x\leftrightarrow y$, and therefore either both are D-finite or neither
is.
Similarly, if one model can be obtained from another by multiplying all multiplicities by a
nonzero constant~$\lambda$, then its generating function can be obtained from the generating
function of the other by sending $t$ to~$\lambda t$, and therefore again either both are D-finite 
or neither is. 

Applying all these filters to the $4^8=65536$ models with possible
multiplicities $a_{u,v}\in\{0,1,2,3\}$ leaves us with 30307 nontrivial models
(including, for the sake of completeness, the 79 interesting models with $a_{u,v}\in\{0,1\}$
that have already been completely classified).

\section{The Group of the Model}

For a fixed model, i.e., for a fixed choice of multiplicities $a_{u,v}\in\set C$, consider the functional
equation~\eqref{eq:funeq}. The group associated to the model acts on this equation. Its elements map
the variables $x$ and $y$ to certain rational functions in $x$ and~$y$, which are chosen in such a way
that all the group elements leave the \emph{kernel polynomial}
\[
  K(x,y,t) := 1 - t \sum_{u,v} a_{u,v}x^u y^v
\]
fixed. It is easy to check that the two particular transformations $\Phi,\Psi\colon\set C(x,y)\to\set C(x,y)$
defined by 
\begin{alignat*}1
  \Phi\colon(x,y)\mapsto \Bigl(\frac1x \frac{\sum_v a_{-1,v}y^v}{\sum_v a_{1,v}y^v},\ y\Bigr),
 \qquad
 &\Psi\colon(x,y)\mapsto \Bigl(x,\ \frac1y \frac{\sum_u a_{u,-1}x^u}{\sum_u a_{u,1}x^u}\Bigr)
\end{alignat*}
have this property. It is also easy to check that $\Phi$ and $\Psi$ are
involutions, i.e., $\Phi^2=\Psi^2=\mathrm{id}$.

The group $G$ is defined as the group generated by $\Phi$ and~$\Psi$ under
composition.

Note that we do not need to worry that one of the denominators $\sum_u
a_{u,1}x^u$ or $\sum_v a_{1,v}y^v$ is identically zero, because this only
happens for models that are uninteresting in the sense of the previous
section. For the same reason, we may also assume that the numerators $\sum_u a_{u,-1}x^u$ and
$\sum_v a_{-1,v}y^v$, respectively, are nonzero polynomials.  In order
to argue that the composition of rational functions into the power series of
equation~\eqref{eq:funeq} is algebraically meaningful, recall that the series in
question belong to $\set Q[x,y][[t]]$, so the result of the composition can be
naturally interpreted as an element of~$\set C(x,y)[[t]]$.

The group $G$ is finite if and only if $(\Phi\Psi)^n=\mathrm{id}$ for some $n\in\set N$,
and this is the case if and only if 
\[
  \begin{matrix}
    G=\bigl\{\!\! & \mathrm{id}, & \Phi\Psi, & (\Phi\Psi)^2, & \cdots, & (\Phi\Psi)^{n-1},\\[2pt]
              & \Phi, & (\Phi\Psi)\Phi, & (\Phi\Psi)^2\Phi, & \cdots, & (\Phi\Psi)^{n-1}\Phi&\!\!\bigr\},
  \end{matrix}
\]
where all the listed elements are distinct. In other words, $G$ is either the dihedral group with
$2n$ elements, or infinite. The sign $\sgn(g)$ of an element $g\in G$ is
defined to be $1$ if $g=(\Phi\Psi)^k$ for some $k$, and $-1$ otherwise. 

\section{Models with Group D4}\label{sec:d4}

As there is obviously no way to choose $a_{u,v}$ such that $\Phi\Psi=(\frac1x
r(y),\frac1y s(x))=(x,y)=\mathrm{id}$, the smallest possible $n\in\set N$ with
$(\Phi\Psi)^n=\mathrm{id}$ is~$2$. The group with $(\Phi\Psi)^2=\mathrm{id}$ is
the dihedral group D4 with four elements.  In order to determine the models
which lead to this group, regard the $a_{u,v}$ as variables and compute
$(p,q):=\Phi(\Psi(x, y)) - \Psi(\Phi(x,y))$. This is a pair of rational functions in
$x,y$ whose coefficients are rational functions in the $a_{u,v}$ over the
rational numbers. Write $p,q$ as quotients of polynomials in $x,y$ whose
coefficients are polynomials in $a_{u,v}$ with integer coefficients. We want to
know the possible choices of $a_{u,v}$ for which $p$ and $q$ become zero. (Note that
$\Phi\Psi=\Psi\Phi\iff (\Phi\Psi)^2=\mathrm{id}$ because $\Phi$ and $\Psi$ are
involutions.) In order to find these~$a_{u,v}$, consider the ideal in $\set
Q[a_{-1,-1},\dots,a_{1,1}]$ generated by the coefficients of all monomials $x^i
y^j$ in the numerator of~$p$ and the coefficients of all monomials $x^i y^j$ in the
numerator of~$q$. This ideal basis consists of 36 homogeneous
polynomials of degree~4, which we don't reproduce here because of its
length. Using Gr\"obner basis techniques~\citep{becker93}, we can determine the
irreducible components of the radical of this ideal. We have used the commands
\verb|facstd| and \verb|minAssGTZ| of the software package Singular~\citep{greuel02}
for this step. It turns out that the two irreducible components are generated by
\begin{alignat*}1
 &\{\,a_{0,1}a_{1,-1} - a_{0,-1}a_{1,1}, \
    a_{-1,1}a_{1,-1} - a_{-1,-1}a_{1,1}, \
    a_{-1,1}a_{0,-1} - a_{-1,-1}a_{0,1}\,\},\text{ and}\\
 &\{\,a_{1,0}a_{-1,1} - a_{-1,0}a_{1,1}, \
    a_{1,-1}a_{-1,1} - a_{-1,-1}a_{1,1}, \
    a_{1,-1}a_{-1,0} - a_{-1,-1}a_{1,0}\,\}.
\end{alignat*}
As the latter is obtained from the former by replacing all $a_{u,v}$ by $a_{v,u}$, it
suffices to consider one of the two components, say the first. The equations in this
component are equivalent to saying that the vectors $(a_{-1,-1},a_{0,-1},a_{1,-1})$ and 
$(a_{-1,1},a_{0,1},a_{1,1})$ are linearly dependent. Since the models where one or 
both of these vectors are zero are uninteresting, the interesting models
leading to the group D4 are precisely those for which there exists a constant
$\lambda\neq0$ such that $a_{-1,v}=\lambda a_{1,v}$ for $v=-1,0,1$. We then have
\[
  \Phi(x, y) = \Bigl(\frac\lambda x, y\Bigr)
  \quad\text{and}\quad
  \Psi(x, y) = \Bigl(x, \frac1y\frac{\lambda a_{1,-1}x^{-1} + a_{0,-1} + a_{1,-1}x}{\lambda a_{1,1}x^{-1} + a_{0,1} + a_{1,1}x}\Bigr) .
\]
At this point, we can proceed analogously to Bousquet-M\'elou and Mishna (cf.\ their Proposition~5): 
multiplying \eqref{eq:funeq} on both sides by~$xy/K(x,y,t)$ and forming the \emph{orbit sum} gives
the general relation
\[
  \sum_{g\in G} \sgn(g)\,g(xy\,f(x,y,t)) = \frac1{K(x,y,t)}\sum_{g\in G} \sgn(g)\,g(xy),
\]
which holds whenever the group is finite. For the special case under consideration, 
the right hand side evaluates to 
\begin{alignat*}1
 &\frac1{K(x,y,t)}\Bigl(
  xy - \frac{y\lambda}x - \frac xy \frac{\lambda a_{1,-1}x^{-1} + a_{0,-1} + a_{1,-1}x}{\lambda a_{1,1}x^{-1} + a_{0,1} + a_{1,1}x} + \frac\lambda{xy} \frac{\lambda a_{1,-1}x^{-1} + a_{0,-1} + a_{1,-1}x}{\lambda a_{1,1}x^{-1} + a_{0,1} + a_{1,1}x}\Bigr)\\
 &= \frac{(x^2-\lambda) (a_{0,1}x y^2 - a_{0,-1}x -(\lambda +x^2) (a_{1,-1}-a_{1,1}y^2))}{x y (a_{1,1}(\lambda +x^2)+ a_{0,1}x) K(x,y,t)}.
\end{alignat*}
For the left hand side, we have
\begin{alignat*}1
  x y \, f(x,y,t)
- \frac{\lambda y}x\, f\Bigl(\frac\lambda x, y, t\Bigr)
- \frac xys(x)\, f\Bigl(x, \frac1 y s(x), t\Bigr)
+ \frac\lambda{xy}s(x)\, f\Bigl(\frac\lambda x, \frac1y s(x), t\Bigr),
\end{alignat*}
where we abbreviate $s(x)=\frac{\lambda a_{1,-1}x^{-1} + a_{0,-1} + a_{1,-1}x}{\lambda a_{1,1}x^{-1} + a_{0,1} + a_{1,1}x}$. The identity holds in $\set Q(x,y)[[t]]$, but it can be seen that all quantities actually belong to
$\set Q(x)[y,y^{-1}][[t]]$. The last two terms of the equation involve only negative exponents with respect to~$y$,
so taking the positive part $[y^>]$ will kill them. The remaining terms happen to belong to $\set Q[x,x^{-1}][[t]]$,
and since the second term only has negative exponents with respect to~$x$, taking the positive part $[x^>]$ will eliminate
it and only leave the first. It follows that
\[
 f(x,y,t) = \frac1{xy}[x^>][y^>] 
\frac{(x^2-\lambda) (a_{0,1}x y^2 - a_{0,-1}x -(\lambda +x^2) (a_{1,-1}-a_{1,1}y^2))}{x y (a_{1,1}(\lambda +x^2)+ a_{0,1}x) K(x,y,t)}.
\]
Alternatively, we could interpret the elements of $\set Q(x,y)[[t]]$ as
elements of multivariate formal Laurent series field $\set Q_{\leq}((x,y,t))$ for a term order
$\leq$ with $x,y\leq 1\leq t$ and do the positive part extraction with respect
to $x$ and $y$ simultaneously. See \cite{aparicio12} for a discussion of formal Laurent
series in several variables. In any case, we can summarize the result of this section as follows.

\begin{thm}
  The interesting quarter plane models whose group is D4 are precisely those where 
  $a_{-1,v}=\lambda a_{1,v}$ for $v=-1,0,1$ and some $\lambda\neq0$. All these models are D-finite.
\end{thm}

\noindent 
\begin{tabular}{@{}|p{.971\hsize}|@{}}\hline
  \family{0}{$a_{0,1}a_{1,-1} = a_{0,-1}a_{1,1}$,\hfill\break
  $a_{-1,1}a_{1,-1} = a_{-1,-1}a_{1,1}$,\hfill\break
  $a_{-1,1}a_{0,-1} = a_{-1,-1}a_{0,1}$}{%
  \begin{tikzpicture}[scale=.3]
    \draw[->] (0,0) -- (-1,-1) node[left] {$\scriptstyle 2$};
    \draw[->] (0,0) -- (-1,0) node[left] {$\scriptstyle -3$};
    \draw[->] (0,0) -- (-1,1) node[left] {$\scriptstyle 5$};
    \draw[->] (0,0) -- (0,-1) node[below] {$\scriptstyle -7$};
    \draw[->] (0,0) -- (0,1) node[above] {$\scriptstyle 13$};
    \draw[->] (0,0) -- (1,-1) node[right] {$\scriptstyle 6$};
    \draw[->] (0,0) -- (1,0) node[right] {$\scriptstyle -9$};
    \draw[->] (0,0) -- (1,1) node[right] {$\scriptstyle 15$};
  \end{tikzpicture}}%
  \\\hline
\end{tabular}

\section{Models with Group D6}

We now determine all the choices for $a_{u,v}$ such that $(\Phi\Psi)^3=\mathrm{id}$. As before,
we compute $(p,q):=\Psi(\Phi(\Psi(x,y)))-\Phi(\Psi(\Phi(x,y)))$ and consider the ideal generated
by the coefficients of the numerators with respect to~$x,y$. The basis consists
of 210 homogeneous polynomials of degree~9. The ideal has 34 irreducible components, 18 of which turn
out to contain only uninteresting models. Of the remaining 16 components, 6 can be discarded
because their solution sets are properly contained in the solution set of others. Of the remaining
10 components, 4 can be discarded because they are reflections of others. This leaves us with 
the following 6 families:

\par\medskip\noindent\kern-.4pt
\begin{tabular}{@{}|p{.472\hsize}|p{.472\hsize}|@{}}\hline
  \family{1a}{$a_{1,1}=a_{-1,-1}=0$,\hfill\break
  $a_{-1,1}a_{1,-1}=a_{-1,0}a_{1,0}=a_{0,1}a_{0,-1}$}
   {\begin{tikzpicture}[scale=.3]
    \draw[->] (0,0) -- (0,-1) node[below] {$\scriptstyle 1/2$};
    \draw[->] (0,0) -- (1,-1) node[right] {$\scriptstyle 1/3$};
    \draw[->] (0,0) -- (1,0) node[right] {$\scriptstyle 1/5$};
    \draw[->] (0,0) -- (0,1) node[above] {$\scriptstyle 2$};
    \draw[->] (0,0) -- (-1,1) node[left] {$\scriptstyle 3$};
    \draw[->] (0,0) -- (-1,0) node[left] {$\scriptstyle 5$};
  \end{tikzpicture}\kern-1em\null}%
  &
  \family{1b}{$a_{1,-1}=a_{-1,1}=0$,\hfill\break
    $a_{-1,0}a_{1,0}=a_{-1,-1}a_{1,1}=a_{0,-1}a_{0,1}$}
  {\begin{tikzpicture}[scale=.3]
      \draw[->] (0,0) -- (0,1) node[above] {$\scriptstyle 1/2$};
      \draw[->] (0,0) -- (1,1) node[right] {$\scriptstyle 1/3$};
      \draw[->] (0,0) -- (1,0) node[right] {$\scriptstyle 1/5$};
      \draw[->] (0,0) -- (0,-1) node[below] {$\scriptstyle 2$};
      \draw[->] (0,0) -- (-1,-1) node[left] {$\scriptstyle 3$};
      \draw[->] (0,0) -- (-1,0) node[left] {$\scriptstyle 5$};
    \end{tikzpicture}\kern-1em\null}%
  \\\hline
  \family{2a}{$a_{1,0}=a_{1,1}=0$,\hfill\break
    $a_{0,-1}a_{-1,1}=2a_{0,1}a_{-1,-1}$, \hfill\break
    $a_{0,-1}^2=4a_{1,-1}a_{-1,-1}$,\hfill\break
    $a_{0,-1}a_{0,1}=2a_{-1,1}a_{1,-1}$}
    {\begin{tikzpicture}[scale=.3]
        \draw[->] (0,0) -- (-1,1) node[above] {$\scriptstyle 7$};
        \draw[->] (0,0) -- (0,1) node[above] {$\scriptstyle 7$};
        \draw[->] (0,0) -- (-1,0) node[left] {$\scriptstyle 5$};
        \draw[->] (0,0) -- (-1,-1) node[below] {$\scriptstyle 1$};
        \draw[->] (0,0) -- (0,-1) node[below] {$\scriptstyle 2$};
        \draw[->] (0,0) -- (1,-1) node[below] {$\scriptstyle 1$};
      \end{tikzpicture}}%
  &
  \family{2b}{$a_{1,0}=a_{1,-1}=0$,\hfill\break
    $a_{0,1}a_{-1,-1}=2a_{0,-1}a_{-1,1}$, \hfill\break
    $a_{0,1}^2 = 4a_{1,1}a_{-1,1}$,\hfill\break
    $a_{0,1}a_{0,-1}=2a_{-1,-1}a_{1,1}$}
  {\begin{tikzpicture}[scale=.3]
      \draw[->] (0,0) -- (-1,-1) node[below] {$\scriptstyle 7$};
      \draw[->] (0,0) -- (0,-1) node[below] {$\scriptstyle 7$};
      \draw[->] (0,0) -- (-1,0) node[left] {$\scriptstyle 5$};
      \draw[->] (0,0) -- (-1,1) node[above] {$\scriptstyle 1$};
      \draw[->] (0,0) -- (0,1) node[above] {$\scriptstyle 2$};
      \draw[->] (0,0) -- (1,1) node[above] {$\scriptstyle 1$};
    \end{tikzpicture}}%
  \\\hline
  \family{3a}{$a_{-1,0}=a_{-1,-1}=0$,\hfill\break
    $a_{0,1}a_{1,-1}=2a_{0,-1}a_{1,1}$, \hfill\break
    $a_{0,1}^2 = 4a_{-1,1}a_{1,1}$,\hfill\break
    $a_{0,1}a_{0,-1}=2a_{1,-1}a_{-1,1}$}
   {\begin{tikzpicture}[scale=.3]
       \draw[->] (0,0) -- (1,-1) node[below] {$\scriptstyle 7$};
       \draw[->] (0,0) -- (0,-1) node[below] {$\scriptstyle 7$};
       \draw[->] (0,0) -- (1,0) node[right] {$\scriptstyle 5$};
       \draw[->] (0,0) -- (1,1) node[above] {$\scriptstyle 1$};
       \draw[->] (0,0) -- (0,1) node[above] {$\scriptstyle 2$};
       \draw[->] (0,0) -- (-1,1) node[above] {$\scriptstyle 1$};
     \end{tikzpicture}}%
  &
  \family{3b}{$a_{-1,0}=a_{-1,1}=0$,\hfill\break
    $a_{0,-1}a_{1,1}=2a_{0,1}a_{1,-1}$, \hfill\break
    $a_{0,-1}^2=4a_{-1,-1}a_{1,-1}$,\hfill\break
    $a_{0,-1}a_{0,1}=2a_{1,1}a_{-1,-1}$}
  {\begin{tikzpicture}[scale=.3]
      \draw[->] (0,0) -- (1,1) node[above] {$\scriptstyle 7$};
      \draw[->] (0,0) -- (0,1) node[above] {$\scriptstyle 7$};
      \draw[->] (0,0) -- (1,0) node[right] {$\scriptstyle 5$};
      \draw[->] (0,0) -- (1,-1) node[below] {$\scriptstyle 1$};
      \draw[->] (0,0) -- (0,-1) node[below] {$\scriptstyle 2$};
      \draw[->] (0,0) -- (-1,-1) node[below] {$\scriptstyle 1$};
    \end{tikzpicture}}%
  \\\hline
\end{tabular}\kern-.4pt\null\par\medskip

Note that the families on the right can be obtained from those on the left by reflection about the horizontal 
axis and the families in the third row can be obtained from those in the second row by reversing all arrows.
The families in the first row are closed under reversing arrows. 

\begin{thm}
  The interesting quarter plane models whose group is D6 are precisely those that 
  belong to one or more of the families described in the table above. 
  All these models are D-finite. 
\end{thm}

The remainder of this section is devoted to the D-finiteness claim of this theorem. 

\subsection{Families 1a, 2a, 3a}

These families can be handled very much like the family in Section~\ref{sec:d4} above. Without going into
further details, we just report the resulting formulas for the generating functions. 

For family 1a, let $\lambda=a_{-1,1}a_{1,-1}=a_{-1,0}a_{1,0}=a_{0,1}a_{0,-1}$. 
If $\lambda\neq0$, then all the $a_{u,v}$ are nonzero (except $a_{1,1}$ and $a_{-1,-1}$ of course). 
In this case, the resulting formula for the generating function is 
  \[
    f(x,y,t) = \frac1{xy}[x^>y^>] 
      \frac{(a_{-1,1}-a_{0,-1}xy^{-2})(a_{1,-1}-a_{-1,0}yx^{-2})(\lambda x y - a_{-1,0}a_{0,-1})}
          {\lambda^2  K(x,y,t) }. 
  \]
Otherwise, if $\lambda=0$ and $a_{-1,1}=0$, then $a_{-1,0}\neq0$ and $a_{0,1}\neq0$ 
(otherwise the model is not interesting), but then
$a_{1,0}=0$ and $a_{0,-1}=0$ (by the defining equations), and then $a_{1,-1}\neq0$ (otherwise again the model is not interesting).
In this case, the resulting formula for the generating function is
  \[
    f(x,y,t) = \frac1{xy}[x^>y^>] 
    \frac{(a_{0,1} -a_{1,-1}xy^{-2} ) (a_{1,-1}- a_{-1,0}x^{-2}y )(a_{0,1}x y-a_{-1,0})}
          {a_{0,1}^2 a^{\vphantom2}_{1,-1} K(x,y,t)}. 
  \]
Finally, if $\lambda=0$ and $a_{-1,1}\neq0$, then $a_{1,-1}=0$ and the only interesting cases have
$a_{1,0}\neq0$, $a_{0,-1}\neq0$, and $a_{-1,0}=a_{0,1}=0$. This case is symmetric to the previous case
and therefore not interesting.

For family 2a, we may assume that $a_{1,-1}\neq0$, because otherwise the model is uninteresting. 
Then we can also assume $a_{0,1}\neq0$, because if $a_{0,1}=0$, then the last defining equation 
would imply $a_{-1,1}=0$, which together with $a_{1,1}=0$ would also render the model not interesting. 
Under the assumption $a_{1,-1}\neq0$, $a_{0,1}\neq0$, the generating function can be expressed as
\[
  f(x,y,t) = \frac1{xy}[x^>y^>] 
\frac{
    P(x,y)
    \bigl(a_{-1,-1}- a_{1,-1}x^2+ a_{-1,1}y^2+ a_{-1,0}y\bigr)
    \bigl(2 a_{0,1}y^2 - 2 a_{1,-1}x - a_{0,-1}\bigr) 
  }{4 x^2 y^3 a_{0,1}^2 a^{\vphantom2}_{1,-1} K(x,y,t)} 
\]
where $P(x,y)=2 a_{-1,-1}-2 a_{0,1} x y^2+ a_{0,-1}x+2 a_{-1,1}y^2+2 a_{-1,0}y$. 

For family 3a, models are interesting only when $a_{-1,1}\neq0$ and $a_{0,-1}\neq0$
and $(a_{1,-1},a_{1,0},a_{1,1})\neq(0,0,0)$. Under these assumptions, we obtain the following
expression for the generating function: 
\[
  f(x,y,t) = \frac1{xy}[x^>y^>] 
\frac{Q(x,y)
  \bigl(a_{0,1}y^2 {-} 2 a_{0,-1} {+} 2 a_{-1,1}y^2\kern-1pt/x\bigr) 
  \bigl(a_{1,1}x^2 y {+} a_{1,0}x^2 {+} a_{1,-1}x^2\kern-1pt/y {-} a_{-1,1}y\bigr)}
  {\bigl(a_{1,-1}+ a_{1,0}y+ a_{1,1}y^2\bigr)
   \bigl(4 a_{-1,1} a_{0,-1} + (2a_{-1,1} + a_{0,1}x\bigr) Q(x,y)) K(x,y,t)} 
\]
where $Q(x,y)=2 a_{1,1} x y^2+2 a_{1,0} x y+2 a_{1,-1} x+ a_{0,1}y^2-2
a_{0,-1}$. Note that in this case the denominator contains nontrivial factors
involving both $x$ and~$y$, so the ad-hoc reasoning used in
Section~\ref{sec:d4}, which also works for the families 1a and~2a, does not work here. 
However, there is no problem if we take the viewpoint of multivariate
Laurent series~\citep{aparicio12}, because all that is needed for the argument to
go through is the property that there exists a term order~$\leq$ so that for all $g\in
G\setminus\{\mathrm{id}\}$ and all positive integers $i,j$ the expansion of
$g(x)^i g(y)^j\in\set C(x,y)$ in the multivariate Laurent series field $\set
C_\leq((x,y))$ contains no terms $x^k y^\ell$ where both $k$ and~$\ell$ are positive.
This turns out to be the case. 

\subsection{Family 1b}\label{sec:kreweras}
For the family~1b there are three cases to distinguish. First, when
$a_{-1,-1}=a_{1,0}=a_{0,1}=0$, then $a_{1,1},a_{-1,0},a_{0,-1}$ all must be
nonzero in order for the model to be interesting.  In this case, the generating
function is 
\[
 f(x,y,t)=k\Bigl(\sqrt[3]{a_{0,-1}a_{1,1}/a_{-1,0}^2}\,x, \sqrt[3]{a_{-1,0}a_{1,1}/a_{0,-1}^2}\,y, \sqrt[3]{\vphantom{/a_{0,-1}^2}a_{-1,0}a_{0,-1}a_{1,1}}\,t\Bigr),
\]
where $k(x,y,t)$ is the generating function for classical Kreweras walks (i.e.,
$a_{1,1}=a_{-1,0}=a_{0,-1}=1$), which is well-known to be algebraic~\citep{kreweras65,bousquet05}.  Secondly,
when $a_{1,1}=a_{-1,0}=a_{0,-1}=0$, algebraicity of the generating function can
be established by a similar argument.  The third case is when
$a_{1,1}$,~$a_{-1,-1}$, $a_{1,0}$, $a_{-1,0}$, $a_{0,1}$,~$a_{0,-1}$ are all
nonzero.  In this case it is impossible to express the generating function in
terms of the generating function for the corresponding model without
multiplicities, known as the double Kreweras model. However, if we let
$f_\lambda(x,y,t)$ be the generating function for the family where
$a_{-1,-1}=a_{-1,0}=a_{0,-1}=1$ and $a_{1,1}=a_{1,0}=a_{0,1}=\lambda\neq0$, then
\[
 f(x,y,t)=
 f_{a_{0,1}a_{1,0}^2/(a_{0,-1}a_{1,1}^2)}\Bigl(\frac{a_{1,1}}{a_{0,1}}x, \frac{a_{1,1}}{a_{1,0}} y,
 \frac{a_{-1,-1}a_{1,1}^2}{a_{0,1}a_{1,0}}t\Bigr)
\]
is the generating function of an arbitrary model of family~1b with $a_{1,1}a_{-1,-1}\neq0$.
It therefore suffices to show that $f_\lambda(x,y,t)$ is D-finite. We will show that it is in 
fact algebraic, following the treatment in Section~6.3 of \cite{bousquet10} step by step with the added parameter~$\lambda$. The orbit sum argument does not work here because the orbit sum turns out to be zero. We therefore sum \eqref{eq:funeq} over only half the orbit to obtain a nonzero expression on both sides. This new expression will be more complicated than in the orbit sum case: in general, it will involve the unknown series $f_\lambda(x,y,t)$, $f_\lambda(x,0,t)$ and $f_\lambda(0,0,t)$. However, by careful coefficient extraction, the algebraicity result is still attainable. 

Writing $A_{v}=\sum_{u}a_{u,v}x^{u}$ for $v=-1,0,1$, the half-orbit sum equation reads
\[
 xyf_{\lambda}(x,y,t)-\frac{1}{\lambda x}f_{\lambda}\Bigl(\frac{1}{\lambda xy},y\Bigr)+\frac{1}{\lambda y}f_{\lambda}\Bigl(\frac{1}{\lambda xy},x\Bigr)=\frac{xy-\frac{1}{\lambda x}+\frac{1}{\lambda y}-2txA_{-1} f_{\lambda}(x,0,t)+t f_{\lambda}(0,0,t)}{K(x,y,t)}.
\]
Next we extract the coefficient of $y^{0}$. In order to do so, we use Lemma~7 from \cite{bousquet10}. That is, we solve $K(x,y,t)=0$ for $y$ in terms of $x$ and~$t$:
writing 
$\Delta(x):=t^2\,x^{-2} - 2(t+2\lambda t^2)x^{-1} + (1-6\lambda t^2) - 2\lambda t(1+2t)x + \lambda^2 t^2 x^2$
for the discriminant of~$K(x,y,t)$, the two solutions are 
$Y_0=\bigl(1-t A_0 -\sqrt{\Delta(x)}\bigr)/(2tA_0)$ and $Y_1=1/(\lambda x Y_0)$. 
The coefficient of $y^n$ in $1/K(x,y,t)$ can be expressed in terms of $Y_0$, $Y_1$, and $\Delta(x)$ via
\[
 [y^n]\frac1{K(x,y,t)}=\frac1{\sqrt{\Delta(x)}}\times\left\{
   \begin{array}{ll}
     Y_0^{-n}& \text{if $n\leq0$}\\
     Y_1^{-n}& \text{if $n\geq0$}
   \end{array}
   \right..
\]
Using these facts, extracting the coefficient of $y^0$ on both sides of the half-orbit sum equation leads to
\begin{equation}\label{eq:xx}
 -\frac{1}{\lambda x}d_{\lambda}\Bigl(\frac{1}{\lambda x},t\Bigr)=\frac{1}{\sqrt{\Delta(x)}}\Bigl(xY_0 -\frac{1}{\lambda x}+\frac{1}{\lambda Y_1}-2tx A_{-1}f_{\lambda}(x,0,t)+tf_{\lambda}(0,0,t)\Bigr),
\end{equation}
where $d_\lambda(x,t):=\sum_{i,n} (f_\lambda)_{i,i,n}x^{i}t^n$ is the generating function for walks ending on the diagonal. 

Now we write $\Delta(x)=\frac{t^2}{Z^2}\Delta_-(x)\Delta_+(x)$, where
\[
 \Delta_+(x)=1 - \frac{2\lambda Z(1 + 2Z + 2\lambda Z^2 + 2\lambda ^2Z^3 + \lambda ^2Z^4)}{(1 -
\lambda Z^2)^2}x + \lambda ^2 Z^2 x^2,\quad 
 \Delta_-(x)=\Delta_+\Bigl(\frac 1 x\Bigr),
\]
and where $Z\in\set Q[\lambda][[t]]$ is defined through $Z = \frac{t(1 + 3\lambda Z^2 + 4\lambda(1+\lambda) Z^3 + 3\lambda ^2Z^4 + \lambda ^3Z^6)}{(1 -\lambda Z^2)^2}$ and $Z(0)=0$. 

Multiplying \eqref{eq:xx} by $A_1\sqrt{\Delta_-(x)}$ and using the explicit expressions for $Y_0$ and $Y_1$ given above, we obtain
\[
 \sqrt{\Delta_-(x)}\Bigl(\frac x t -\frac{1}{\lambda x}A_1 d_{\lambda}\Bigl(\frac{1}{\lambda x},t\Bigr)\Bigr)=
 \frac{ZA_1}{t\sqrt{\Delta_+(x)}} \Bigl(\frac{x\left(1-t A_0\right)}{tA_1}-\frac{1}{x\lambda }-2tA_{-1}f_{\lambda}(x,0,t)+tf_{\lambda}(0,0,t)\Bigr).
\]
From this equation, we extract the coefficient of~$x^0$. 
Using $[x^0]d_\lambda(\frac1{\lambda x},t)=[x^0]f_\lambda(x,0,t)=f_\lambda(0,0,t)$, we find
\[
 f_{\lambda}(0,0,t)=\frac{Z-4\lambda Z^3-2\lambda Z^4-2\lambda ^2Z^4-\lambda ^2Z^5}{t(1-\lambda Z^2)^2}. 
\]
With this knowledge, we can now extract the positive part in $x$ on both sides of the same equation to obtain
\begin{alignat*}1
 &f_{\lambda}(x,0,t)= \frac{ x^2 (\lambda  Z^2-1)+2 x
   Z (\lambda  Z+1)-\lambda  Z^3+Z}{2 \lambda  t x (x+1)^2 Z
   (1-\lambda  Z^2)}\sqrt{\Delta_+(x)}\\
 &{}-\frac{Z}{2t(1+x)}\left(\frac{\lambda  t x^3+2 t x+t-x^2}{\lambda  t x (x+1) Z}+\frac{2 (\lambda ^2 Z^3+\lambda  (Z+3)
   Z^2-1)}{(1-\lambda  Z^2)^2}+1\right). 
\end{alignat*}
Noting that $f_{\lambda}(0,y,t)=f_{\lambda}(y,0,t)$, we conclude from equation~\eqref{eq:funeq} that $f_{\lambda}(x,y,t)$ is algebraic.

\subsection{Families 2b, 3b}\label{sec:new}

The orbit sum argument also fails for these families. For the models in
family~2b the orbit sum is zero, while in family~3b the orbit sum is nonzero but
the desired term $f(x,y,t)$ cannot be isolated by taking the positive part
because there are group elements~$g\neq\mathrm{id}$ for which $f(g(x),g(y),t)$
also contributes terms with positive exponents to the orbit sum. Because of the lack
of symmetry, the half orbit sum argument used for family~1b does not seem to
apply either.

One model from each of these two families were already encountered by~\cite{bostan14},
and computer proofs have been given there that the generating function for the
model belonging to family~2b is algebraic and the model belonging to family~3b
is (transcendental) D-finite. The models considered by~\cite{bostan14} are
$a_{1,0}=a_{1,-1}=a_{-1,0}=0$, $a_{-1,1}=\tfrac12a_{0,1}=a_{1,1}=a_{-1,-1}=a_{0,-1}=1$ (case~2b), and its reverse 
$a_{-1,0}=a_{-1,1}=a_{1,0}=0$, $a_{1,-1}=\tfrac12a_{0,-1}=a_{-1,-1}=a_{1,1}=a_{0,1}=1$ (case~3b).

We were able to extend these computer proofs to the more general cases where
$a_{-1,0}=\lambda$ (case~2b), and $a_{1,0}=\lambda$ (case~3b), respectively, are
formal parameters. From here, every other model of the respective family can be
reached by an appropriate algebraic substitution: if $f_\lambda(x,y,t)$ is the 
generating function for the model $a_{1,0}=a_{1,-1}=a_{-1,0}=0$,
$a_{-1,1}=\tfrac12a_{0,1}=a_{1,1}=a_{-1,-1}=a_{0,-1}=1$, $a_{-1,0}=\lambda$, then
\[
  f(x,y,t) = f_{a_{-1,0}/\sqrt{a_{-1,-1}a_{-1,1}}}\Bigl(\frac{a_{0,-1}}{a_{-1,-1}}\,x,\ \sqrt{\frac{a_{-1,1}}{a_{-1,-1}}}\,y,\ \frac12\sqrt{\frac{a_{-1,-1}}{a_{-1,1}}}\, t\Bigr)
\]
is the generating function for an arbitrary model of family~2b, and likewise for 
family~3b. (Models where the $a_{u,v}$'s appearing in the denominators are zero are not 
interesting.)

The computational techniques we used were introduced by
\cite{kauers07v,kauers09b,bostan10}, and they have been described for the cases
$\lambda=0$ in the paper of~\cite{bostan14}. We do not repeat these
explanations again but only remark that the additional symbolic
parameter~$\lambda$ has made the calculations considerably more expensive. 
The computations were done using software of \cite{kauers09a} and \cite{koutschan10c}. 
The bottleneck was the construction of a certified recurrence for
$(f_\lambda)_{0,0,n}$. The (nonminimal) recurrence we found has order~14 and
degrees 30,~26 in~$n$,~$\lambda$, respectively; the certificate for this recurrence is 16
gigabytes long! From this recurrence it can be deduced that $f_\lambda(0,0,t)$
is the unique formal power series $T\in\set Q[\lambda][[t]]$ with $T(0)=1$ and
\[
 t^4T^2 
+(2 t \lambda+1)t^2T
+t (4 t+1)-(3 t^2 (\lambda-4)+3 t+1)Z+t (6 t+1) (\lambda+2) Z^2=0,
\]
where $Z\in\set Q[\lambda][[t]]$ is the unique formal power series with $Z(0)=0$ and
$t=\frac{Z(4Z+1)}{1+6Z+12Z^2+4(2+\lambda)Z^3}$. Using this equation and the
functional equation~\eqref{eq:funeq} (with $f_\lambda$ in place of~$f$), we
could then prove the correctness of guessed polynomial equations
$P(x,t,\lambda,f_\lambda(x,0,t))=Q(y,t,\lambda,f_\lambda(0,y,t))=0$, which in
turn can be used to deduce that $f_\lambda(x,0,t)$ is the unique formal power
series $U\in\set Q[x,\lambda][[t]]$ with $U(0)=1$ and
\[
 (x+1)^2t^4 U^2
+ (2 t \lambda-x+1)t^2U
+t (t (x+4)+1)-(3 t^2 (\lambda-4)+3 t+1)Z+t (6 t+1) (\lambda+2) Z^2=0
\]
and that $f_\lambda(0,y,t)$ is $\bigl(-1+\sqrt{\mathstrut 1+tyV}\bigr)/(ty)$ where $V$ 
the unique formal power series $V\in\set Q[y,\lambda][[t]]$ with
$V(0)=1$ and 
\begin{alignat*}1
  &(\lambda y+y^2+1)^2 t^4 V^2+\bigl(4 t^2 (6 t+1) (\lambda+2)y Z^2 -4 t y (3 t^2 \lambda-12 t^2+3 t+1)Z\\
  &\quad{}+t (6 t^2 \lambda y^2+4 t^2 \lambda+2 t^2 y^3+18 t^2 y-2 t \lambda y+2 t y^2+4 t y+2 t-y)\bigr)V\\
  &{}+t (4 t \lambda y+t y^2+16 t+2 y+4)-4 (3 t^2 \lambda-12 t^2+3 t+1)Z+4 (6 t+1) t (\lambda+2) Z^2
     = 0.
\end{alignat*}
Together with the functional equation, it finally follows that
$f_\lambda(x,y,t)$ is algebraic.

For the generating function $\bar f_\lambda(x,y,t)$ of the model with
$a_{-1,0}=a_{-1,1}=a_{0,1}=0$,
$a_{-1,-1}=\tfrac12a_{0,-1}=a_{1,-1}=a_{1,1}=a_{0,1}=1$, $a_{1,0}=\lambda$ from
family~3b, we have that $\bar f_\lambda(0,0,t)=f_\lambda(0,0,t)$ (for
combinatorial reasons), and we can use this and the functional equation to
certify guessed systems of partial linear differential equations for $\bar
f_\lambda(x,0,t)$ and $\bar f_\lambda(0,y,t)$ which then together with the
function equation~\eqref{eq:funeq} (now with $\bar f_\lambda$ in place of~$f$)
imply that $\bar f_\lambda(x,y,t)$ is D-finite.  The equations are somewhat too
large to be included here: $\bar f_\lambda(x,0,t)$ satisfies a differential
equation of order~11 with respect to~$t$ with polynomial coefficients of respective degrees
82,~90,~110 in $x$,~$\lambda$ and~$t$, while $\bar f_\lambda(0,y,t)$
satisfies a differential equation of order~11 with respect to~$t$ with
polynomial coefficients of respective degrees 70,~58,~90 in $y$,~$\lambda$ and~$t$.

\section{Models with Group D8}\label{sec:d8}

For the possible values of $a_{u,v}$ such that $(\Phi\Psi)^4=\mathrm{id}$, we obtain, after
discarding components that only contain uninteresting models or are redundant or are reflections
of others, three essentially different prime ideals. One of them is the ideal from Section~\ref{sec:d4},
which appears again because $(\Phi\Psi)^2=\mathrm{id}$ implies $(\Phi\Psi)^4=\mathrm{id}$. The
other two define the following families:

\par\medskip\noindent\kern-.4pt
\begin{tabular}{@{}|p{.472\hsize}|p{.472\hsize}|@{}}\hline
  \family{4a}{$a_{1,-1}a_{-1,1}=a_{1,0}a_{-1,0}$,\hfill\break
    $a_{1,1}=a_{0,1}=a_{0,-1}=a_{-1,-1}=0$}{%
    \begin{tikzpicture}[scale=.3]
      \draw[->] (0,0) -- (-1,0) node[left] {$\scriptstyle 6$};
      \draw[->] (0,0) -- (-1,1) node[left] {$\scriptstyle 3$};
      \draw[->] (0,0) -- (1,0) node[right] {$\scriptstyle 2$};
      \draw[->] (0,0) -- (1,-1) node[right] {$\scriptstyle 4$};
    \end{tikzpicture}}%
  &
  \family{4b}{$a_{1,1}a_{-1,-1}=a_{1,0}a_{-1,0}$,\hfill\break
    $a_{1,-1}=a_{1,0}=a_{-1,0}=a_{-1,1}=0$}{%
    \begin{tikzpicture}[scale=.3]
      \draw[->] (0,0) -- (1,0) node[right] {$\scriptstyle 2$};
      \draw[->] (0,0) -- (1,1) node[right] {$\scriptstyle 4$};
      \draw[->] (0,0) -- (-1,0) node[left] {$\scriptstyle 6$};
      \draw[->] (0,0) -- (-1,-1) node[left] {$\scriptstyle 3$};
    \end{tikzpicture}}%
  \\\hline
\end{tabular}\kern-.4pt\null\par\medskip

In family~4a, we must have $a_{1,-1}\neq0$ and $a_{-1,1}\neq0$ for a model
to be interesting. But then $a_{1,-1}a_{-1,1}\neq0$ implies also $a_{1,0}\neq0$
and $a_{-1,0}\neq0$ through the first defining equation. Similarly, we can 
assume for the models in family~4b that $a_{1,0},a_{-1,0},a_{1,1}a_{-1,-1}$ all
are nonzero. 

For family~4a, the orbit sum argument applies and yields
\[
  f(x,y,t) = \frac1{xy}[x^>y^>]
  \frac{(a_{1,-1}x/y - a_{-1,1}y/x)
        (a_{1,-1}/y - a_{1,0})
        (a_{1,0}x - a_{-1,0}/x)
        (a_{1,0}x - a_{-1,1}y/x)}
      {a_{-1,1}a_{1,0}^3 K(x,y,t)}  
\]
as expression for the generating function. 

For family~4b the orbit sum is zero, but it was pointed out by \cite{bostan14} in their Section~6.2
that its D-finiteness can be deduced from the D-finiteness of the corresponding model without multiplicities. 
Indeed, if $g(x,y,t)$ denotes the generating function
for the Gessel model $\{\swarrow,\leftarrow,\to,\nearrow\}$ without multiplicities,
then we have
\[
  f(x,y,t) = g\Bigl(\sqrt{\frac{a_{1,0}}{a_{-1,0}}}\,x, \ \frac{a_{1,1}}{a_{1,0}}\,y, \ \sqrt{\mathstrut a_{1,0}a_{-1,0}}\,t\Bigr)
\]
for the general generating function of models of family~4b. Since $g(x,y,t)$ is
known to be algebraic \citep{bostan10,bostan13}, it follows that all the models of family~4b
are algebraic. 

\begin{thm}
  The interesting quarter plane models whose group is D8 are precisely those that 
  belong to one of the families described in the table above. 
  All these models are D-finite. 
\end{thm}

\section{Models with Larger Groups}

For $n\geq5$ we failed to compute the prime decomposition of the ideal of
relations among the $a_{u,v}$ that ensures $(\Phi\Psi)^n=\mathrm{id}$. The
required calculations become too expensive. However, in our search over all the
30307 quarter plane models with multiplicites in $\{0,1,2,3\}$ we did encounter,
very much to our surprise, the following three models that do not belong to any
of the families discussed so far. Their group is~D10.
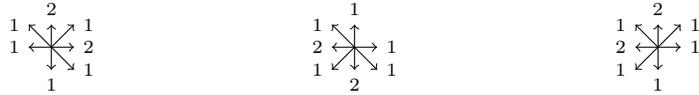
\begin{center}
  \begin{tikzpicture}[scale=.3] 
    \draw[->] (0,0) -- (-1,0) node[left] {$\scriptstyle 1$};
    \draw[->] (0,0) -- (-1,1) node[left] {$\scriptstyle 1$};
    \draw[->] (0,0) -- (0,1) node[above] {$\scriptstyle 2$};
    \draw[->] (0,0) -- (1,1) node[right] {$\scriptstyle 1$};
    \draw[->] (0,0) -- (1,0) node[right] {$\scriptstyle 2$};
    \draw[->] (0,0) -- (1,-1) node[right] {$\scriptstyle 1$};
    \draw[->] (0,0) -- (0,-1) node[below] {$\scriptstyle 1$};
  \end{tikzpicture}\hfil
  \begin{tikzpicture}[scale=.3] 
    \draw[->] (0,0) -- (-1,-1) node[left] {$\scriptstyle 1$};
    \draw[->] (0,0) -- (-1,0) node[left] {$\scriptstyle 2$};
    \draw[->] (0,0) -- (-1,1) node[left] {$\scriptstyle 1$};
    \draw[->] (0,0) -- (0,1) node[above] {$\scriptstyle 1$};
    \draw[->] (0,0) -- (1,0) node[right] {$\scriptstyle 1$};
    \draw[->] (0,0) -- (1,-1) node[right] {$\scriptstyle 1$};
    \draw[->] (0,0) -- (0,-1) node[below] {$\scriptstyle 2$};
  \end{tikzpicture}\hfil
  \begin{tikzpicture}[scale=.3] 
    \draw[->] (0,0) -- (-1,-1) node[left] {$\scriptstyle 1$};
    \draw[->] (0,0) -- (-1,0) node[left] {$\scriptstyle 2$};
    \draw[->] (0,0) -- (-1,1) node[left] {$\scriptstyle 1$};
    \draw[->] (0,0) -- (0,1) node[above] {$\scriptstyle 2$};
    \draw[->] (0,0) -- (1,0) node[right] {$\scriptstyle 1$};
    \draw[->] (0,0) -- (1,1) node[right] {$\scriptstyle 1$};
    \draw[->] (0,0) -- (0,-1) node[below] {$\scriptstyle 1$};
  \end{tikzpicture}
\end{center}
The orbit sum is zero, and guessing suggests that for all three models the generating function is algebraic. 

The models on the left and in the middle can be obtained from one another by reversing 
arrows, therefore these two models have the same number of walks returning to the origin. 
If $Z\in\set Q[[t]]$ is the unique formal power series with $Z(0)=0$ satisfying
$Z = t(4 Z^3+8 Z^2+2 Z+1)$, 
so that $Z=t+2t^2+12t^3+60t^4+\cdots$, then we believe that 
$f(0,0,t) = \frac Zt (1 - 2Z + 2Z^3)$. More generally, for the model
on the left we seem to have
\[
  f(x,0,t)=
  f(0,x,t)=
  \frac{P(x,Z)-(x-2 Z) (2 xZ+x-1) \sqrt{1-4 xZ (Z+1)}}{2 t x^2 (x+1)^2 Z}
\]
with $P(x,Z)=2Z + (x-1) x (4 Z^3+8 Z^2-2(x-1)Z+1)$, and then, using equation~\eqref{eq:funeq}, 
the full generating function $f(x,y,t)$ can be expressed in terms of all these algebraic series. 
For the model in the middle, we find a slightly messier expression for $f(x,0,t)=f(0,x,t)$, also
quadratic over $t,Z,x$, which again implies (if correct) that $f(x,y,t)$ is algebraic.
Since the first two models are symmetric about the diagonal, we expect that
these formulas can be proven in a similar way as the models of family~1b in
Section~\ref{sec:kreweras} above, but we have not gone through the details of the
required calculations.

The model on the right also seems to have an algebraic generating function. We
found that $f(x,0,t)$ seems to satisfy an algebraic equation $P(x,t,f(x,0,t))=0$
for some irreducible polynomial $P\in\set Z[x,t,T]$ of respective degrees 40,
45, 24 in $x,t,T$, and $f(0,y,t)$ seems to satisfy an algebraic equation
$Q(y,t,f(0,y,t))=0$ for some irreducible polynomial $Q\in\set Z[y,t,T]$ of
respective degrees 64, 45, 24 in $y,t,T$. We expect that these algebraic
equations can be proven by computer algebra in a similar way as the models of
family~2b in Section~\ref{sec:new} above, but this would require immense
calculations which we have not carried out.

The model obtained from the model on the right by reversing arrows is just its 
reflection and therefore also algebraic but not of interest. 

Using substitutions like in earlier sections, the three models can be used
to generate three families of models. The corresponding ideals of defining
relations for the $a_{u,v}$ have dimension three. We do not know whether these
families completely characterize all the interesting models whose group is~D10,
nor do we know anything about models for even larger groups. Does there exist
for every $n\geq2$ a quarter plane model with multiplicities whose group is~D$2n$?


\bibliographystyle{abbrvnat}
\bibliography{sample}

\end{document}